\documentclass[12pt]{amsart}
\usepackage{amssymb}
\usepackage{amsmath,amssymb,amsfonts,amsthm,latexsym}
\usepackage{mathrsfs}
\usepackage{color}
\usepackage{hyperref}
\usepackage{graphicx}

\theoremstyle{plain}
\newtheorem{theorem}{Theorem}

\newtheorem{lemma}{Lemma}
\newtheorem{proposition}{Proposition}

\theoremstyle{definition}
\newtheorem{definition}{Definition}

\theoremstyle{example}

\theoremstyle{remark}

\numberwithin{equation}{section}

\begin{document}
\begin{center}
{\bf\Large Combinatorial analysis of interacting RNA molecules}
\\
\vspace{15pt} Thomas J. X. Li, Christian M. Reidys$^{\,\star}$
\end{center}

\begin{center}
Center for Combinatorics, LPMC-TJKLC\\
Nankai University  \\
Tianjin 300071\\
         P.R.~China\\
         Phone: *86-22-2350-6800\\
         Fax:   *86-22-2350-9272\\
duck@santafe.edu
\end{center} 
\centerline{\bf Abstract}{\small
Recently several minimum free energy (MFE) folding algorithms for
predicting the joint structure of two interacting RNA molecules have
been proposed. Their folding targets are interaction structures, that
can be represented as diagrams with two backbones drawn horizontally on
top of each other such that (1) intramolecular and intermolecular
bonds are noncrossing and (2) there is no ``zig-zag'' configuration.
This paper studies joint structures with arc-length at least four in
which both, interior and exterior stack-lengths are at least two (no
isolated arcs). The key idea in this paper is to consider a new type
of shape, based on which joint structures can be derived via
symbolic enumeration. Our results imply simple asymptotic formulas
for the number of joint structures with surprisingly small
exponential growth rates. They are of interest in the context of
designing prediction algorithms for RNA-RNA interactions.}




\section{Introduction}


RNA-RNA binding is an important phenomenon observed in various classes
of non-coding RNAs and plays a crucial role in a number of regulation
processes. Regulatory antisense RNAs control gene expression by prohibiting
the translation of a target mRNA through establishing stable base pairing
interactions. Examples include the regulation of translation in both:
prokaryotes~\cite{Vogel:07} and eukaryotes~\cite{McManus,Banerjee},
the targeting of chemical modifications~\cite{Bachellerie},
insertion editing~\cite{Benne}, and transcriptional control
\cite{Kugel}. More and more evidence suggests, that RNA-RNA interactions also
play a role for the functionality of long mRNA---like ncRNAs~\cite{Hekimoglu}.
A common theme in many RNA classes, including miRNAs,
snRNAs, gRNAs, snoRNAs, and in particular many of the procaryotic small
RNAs, is the formation of RNA-RNA interaction structures that are much
more complex than simple complementary sense-antisense interactions.
The interaction between two RNAs is governed by the same physical
principles that determine RNA folding: the formation of specific
base pairs patterns whose energy is largely determined by base
pair stacking and loop strains.
As a result, secondary structures are an appropriate level of description
to quantitatively understand the thermodynamics of RNA-RNA binding.

\begin{figure}
\begin{center}
\includegraphics[width=0.9\textwidth]{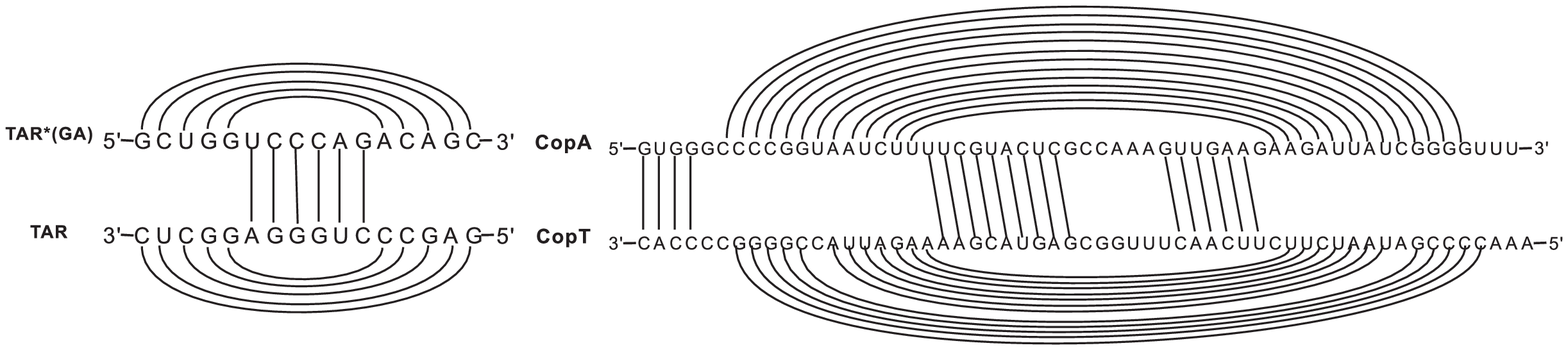}
\end{center}
\caption{Natural joint structures. Known interaction bonds of
$\text{TAR}^{\ast}$(GA)-TAR~\cite{Lebars} and \textsf{CopA}-\textsf{CopT}
\cite{Salari} are displayed.
}\label{F:JSintro}
\end{figure}

Alkan \emph{et al.}~\cite{Alkan:06} proved that the RNA-RNA interaction
prediction (RIP) problem is in its general form NP-complete.
Nevertheless, we are facing increasing demand for efficient
computational methods for RIP.
By restricting the space of allowed configurations, polynomial-time algorithms
on secondary structure level have been derived.
Pervouchine~\cite{Pervouchine:04} and Alkan \emph{et al.}~\cite{Alkan:06}
proposed MFE folding algorithms for predicting the \emph{joint structure} of
two interacting RNA molecules. In this model, ``joint structure'' means
that the intramolecular structures of each partner are
pseudoknot-free, that the intermolecular binding pairs are
noncrossing, and that there is no so-called ``zig-zag''
configuration, see Fig.~\ref{F:JSintro}.
The optimal joint structure can be computed in $O(N^6)$ time and $O(N^4)$ space
by means of dynamic programming~\cite{Alkan:06,Pervouchine:04,rip:09,Backofen}.
Extensions involving the partition function were proposed by Chitsaz \emph{et al.}
(\texttt{piRNA})~\cite{Backofen} and Huang \emph{et al.} (\texttt{rip})
\cite{rip:10}.
In contrast to the situation for RNA secondary structures
\cite{Waterman:78,Waterman:94a}, much less is known about joint
structures. Only joint structures of arc-length greater than or
equal to two have been studied in~\cite{Li:10}. However, the
biochemistry of nucleotide-pairings, favors parallel stacking of
bonds due to entropy and the minimum length of intramolecular bonds
of four. Unfortunately, the biophysically relevant class of
canonical joint structures with arc-length $\geq 4$, is not governed
by the framework in~\cite{Li:10}.

In this paper, we introduce the general framework dealing with
$\sigma$-canonical joint structures having arc-length $\geq
\sigma+2$. In particular, our results apply to the class of
canonical joint structures having arc-length $\ge 4$. Our results
are relevant for the design and analysis of RIP folding algorithms
and show that the numbers of $\sigma$-canonical joint structures
with arc-length $\geq \sigma+2$ exhibit surprisingly small
exponential growth rates.

This paper is organized as follows: in Section~\ref{S:JS} we
introduce joint structures along the lines of~\cite{rip:09} and in
Section~\ref{S:ReShape} we compute, along the lines of
\cite{Modular}, the generating function of refined shapes via
symbolic enumeration. In Section~\ref{S:inflate} we show how to
inflate refined shapes into joint structures and derive the
generating function of joint structures. Section~\ref{S:asym}
presents the singularity analysis and asymptotic formulas. We
finally integrate our results in Section~\ref{S:discussion}.


\section{Secondary structures and joint structures}\label{S:JS}


Let us begin by discussing some basic results of
\cite{Waterman:78,Waterman:94a,Jin:09}. Let $f(n)$ denote the number
of all noncrossing matchings of $n$ arcs having the generating
function ${\bf F}(z)= \sum f(n)\, z^n$. Recursions for $f(n)$ allow
us to derive $z\, {\bf F}(z)^2 -{\bf F}(z) + 1 = 0$, that is we have
\begin{equation*}
{\bf F}(z) = \frac{1-\sqrt{1-4z}}{2z}.
\end{equation*}
Let $\mathcal{T}_{\sigma}^{[\lambda]}$ denote the combinatorial
class of $\sigma$-canonical secondary structures having arc-length
$\geq \lambda$ and $T_{\sigma}^{[\lambda]}(n)$ denote the number of
all $\sigma$-canonical secondary structures over $n$ vertices having
arc-length $\geq \lambda$ and
\begin{equation*}
{\bf T}^{[\lambda]}_{\sigma}(z)= \sum T_{\sigma}^{[\lambda]}(n)\, z^n.
\end{equation*}
\begin{theorem}\label{T:SecondSigLamb}\cite{Jin:09}
Let $\sigma\in\mathbb{N}$, $z$ be an indeterminant and let
\begin{eqnarray*}
u_\sigma(z) & = & \frac{(z^2)^{\sigma-1}}{z^{2\sigma}-z^2+1}, \\
v_\lambda(z) & = & 1-z+u_\sigma(z) {\sum_{h=2}^\lambda z^h},
\end{eqnarray*}
then, ${\bf T}^{[\lambda]}_{\sigma}(z)$, the generating function of
$\sigma$-canonical structures with minimum arc-length $\lambda$ is given by
\begin{equation*}
{\bf T}_{\sigma}^{[\lambda]}(z) = \frac{1}{v_\lambda(z)} {\bf
F}\left(\left(\frac{\sqrt{u_\sigma(z)}\,z}{v_\lambda(z)}
\right)^2\right),
\end{equation*}
where
\begin{equation*}
{\bf F}(z)   = \frac{1-\sqrt{1-4z}}{2z}.
\end{equation*}
Furthermore
\begin{equation*}
T_{\sigma}^{[\lambda]}(n) \sim  c_{\sigma}^{[\lambda]} \, n^{-\frac{3}{2}}
\left( \frac{1}{\zeta_{\sigma}^{[\lambda]}} \right)^n,
\end{equation*}
where $\zeta_{\sigma}^{[\lambda]}$ is the dominant singularity of
${\bf T}_{\sigma}^{[\lambda]}(z)$ and the minimal positive real solution of
the equation
\begin{equation*}
\left(\frac{\sqrt{u_\sigma(z)}\,z}{v_\lambda(z)} \right)^2 = \frac{1}{4}.
\end{equation*}
\end{theorem}
Theorem~\ref{T:SecondSigLamb} implies that for any
specified $\lambda$ and $\sigma$, ${\bf T}^{[\lambda]}_{\sigma}(z)$ is algebraic
over the rational function field $\mathbb{C}(z)$, since ${\bf F}(z)$ is algebraic
and $v_\lambda(z), u_\sigma(z)$ are both rational functions.

\begin{figure}
\begin{center}
\includegraphics[width=0.7\textwidth]{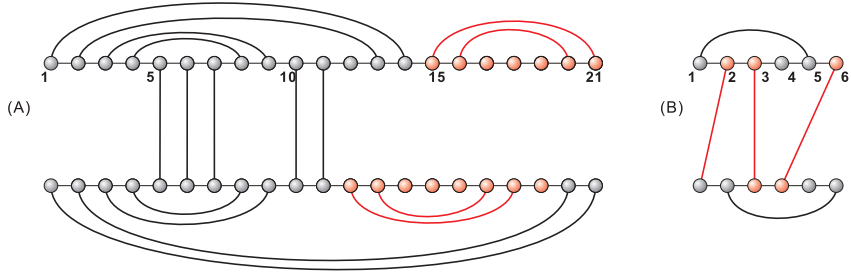}
\caption{\textsf{(A)}: A joint structure $J(R,S,I)$ with arc-length
$\geq 4$ and stack-length $\geq 2$.
Secondary segments (red): the subgraphs $R[15,21]$ and $S[12,19]$.
Ancestors and descendants: for the exterior arc $R_5 S_5$, we have the
following sets of $R$-ancestors and $S$-ancestors of $R_5 S_5$: $\{ R_1 R_{14},
R_2 R_{13}, R_3 R_{9}, R_4 R_{8}\}$ and $\{ S_1 S_{21}, S_2 S_{20}, S_3 S_{9},
S_4 S_{8}\}$. The exterior arc $R_5 S_5$ is a common descendant of $R_1 R_{14}$
and $S_3 S_{9}$, while $R_{10} S_{10}$ is not.
Subsumed arcs: $R_1 R_{14}$ subsumes $S_3 S_9$ and $S_1 S_{21}$.
\textsf{(B)}: A zigzag, generated by $R_2 S_1$, $R_3 S_3$ and $R_6 S_4$.
}
\label{F:JSexample}
\end{center}
\end{figure}

Given two RNA sequences $R= (R_i)_1^{n}$ and $S=(S_j)_1^{m}$ with $n$ and
$m$ vertices, we index the vertices such that $R_1$ is the $5'$ end of $R$ and $S_1$
is the $3'$ end of $S$.  The intramolecular base pair can be represented by an arc
(interior), with its two endpoints contained in either $R$ or $S$. Similarly,
the extramolecular base pair can be represented by an arc (exterior)
with one of its endpoints contained in $R$ and the other in $S$.
When representing arc-configurations, we draw all $R$-arcs in the
upper-halfplane and all $S$-arcs in the lower-halfplane, see
Fig.~\ref{F:JSexample},~{\bf (A)}.

We refer to the subgraph induced by $\{R_i,\ldots,R_j\}$ by $R[i,j]$.
The subgraph $R[i,j]$ ($S[i',j']$) is called \emph{secondary segment} if
there is no exterior arc $R_k S_{k'}$ such that $i \leq k \leq j$
($i' \leq k' \leq j'$), see Fig.~\ref{F:JSexample}, {\bf (A)}. An
interior arc $R_i R_j$ is an $R$-\emph{ancestor} of the exterior arc $R_k
S_{k'}$ if $i<k<j$. Analogously, $S_{i'} S_{j'}$ is an $S$-ancestor
of $R_k S_{k'}$ if $i'<k'<j'$. We also refer to $R_k S_{k'}$ as a
\emph{descendant} of $R_i R_j$ and $S_{i'} S_{j'}$ in this situation, see
Fig.~\ref{F:JSexample}, {\bf (A)}. Furthermore, we call $R_i R_j$
and $S_{i'} S_{j'}$ \emph{dependent} if they have a common descendant and
\emph{independent}, otherwise. Let $R_i R_j$ and $S_{i'} S_{j'}$ be two dependent
interior arcs. Then $R_i R_j$ \emph{subsumes} $S_{i'} S_{j'}$, or $S_{i'} S_{j'}$
is subsumed in $R_i R_j$, if for any $R_k S_{k'}\in
I$, $i'<k'<j'$ implies $i<k<j$, that is, the set of descendants of $S_{i'} S_{j'}$
is contained in the set of descendants of $R_i R_j$,
see Fig.~\ref{F:JSexample},~{\bf (A)}. A \emph{zigzag} is a subgraph
containing two dependent interior arcs $R_{i} R_{j}$ and
$S_{i'} S_{j'}$ neither one subsuming the other, see
Fig.~\ref{F:JSexample}, {\bf (B)}.

A \emph{joint structure}~\cite{Pervouchine:04,Alkan:06,Backofen,rip:09} $J(R,S,I)$, see
Fig.~\ref{F:JSexample},~{\bf (A)}, is a graph such that
\begin{enumerate}
\item $R$, $S$ are secondary structures (each nucleotide being paired with
    at most one other nucleotide via hydrogen bonds, without internal pseudoknots);
\item$I$ is a set of exterior arcs without external pseudoknots, i.e., if
$R_i S_j$, $R_{i'} S_{j'} \in I$ then $i<i'$ implies $j<j'$;
\item $J(R,S,I)$ contains no zig-zags.
\end{enumerate}
We next specify some notations
\begin{itemize}
\item an interior arc (or simply arc) of length $\lambda$ is an arc $R_iR_j$
    ($S_{i'}S_{j'}$) where $j-i=\lambda$ ($j'-i'=\lambda$),
\item an interior stack (or simply stack) of length $\sigma$ is a maximal
      sequence of ``parallel'' interior arcs,
    \begin{eqnarray*}
    (R_i R_j, R_{i+1} R_{j-1}, \ldots , R_{i+\sigma-1} R_{j-\sigma+1} ) &
\text{\rm or}& \\
    (S_i S_j, S_{i+1} S_{j-1}, \ldots , S_{i+\sigma-1} S_{j-\sigma+1} ), &&
    \end{eqnarray*}
\item an exterior stack of length $\tau$ is a maximal sequence of ``parallel''
     exterior arcs,
    \begin{equation*}
    (R_i S_{i'}, R_{i+1} S_{i'+1}, \ldots , R_{i+\tau-1} S_{i'+\tau-1} ).
    \end{equation*}
\end{itemize}
A $\sigma$-\emph{canonical} joint structure is a joint structure with stack-length
$\geq \sigma$. In Fig.~\ref{F:JSexample},
{\bf (A)}, we give an example of $2$-canonical joint
structure with arc-length $\geq 4$.

\begin{figure}
\begin{center}
\includegraphics[width=0.65\textwidth]{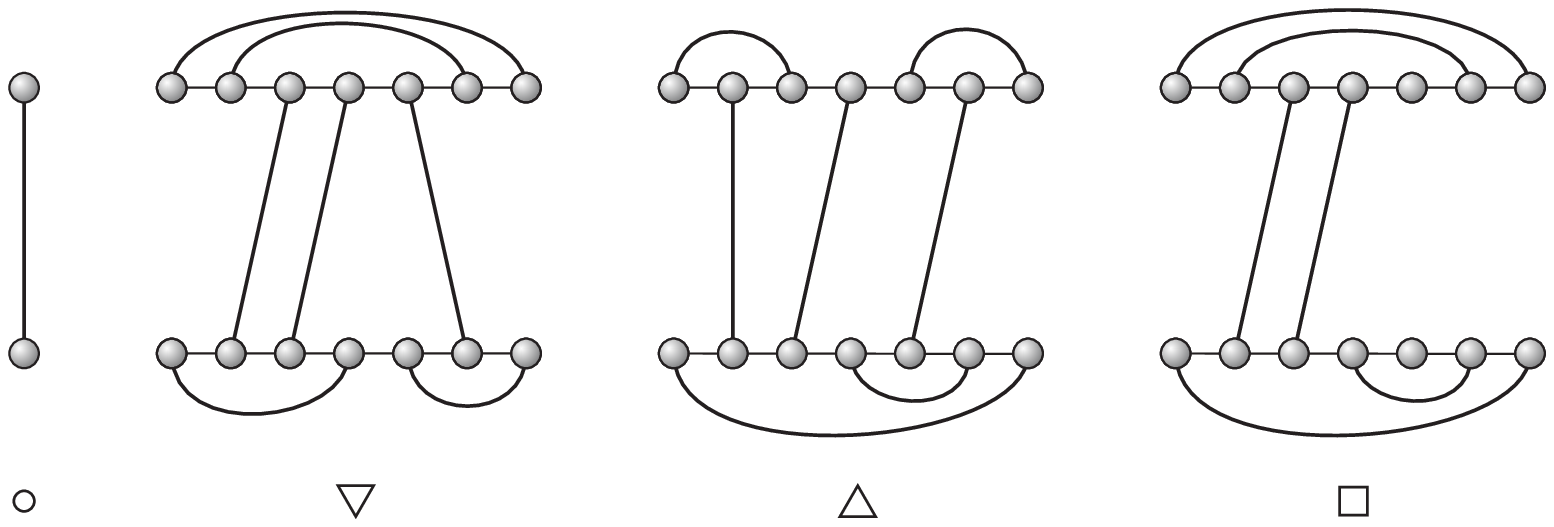}
\caption{The four types of tight structures are defined as follows:
$\circ$ $\colon$ $\{R_i S_{i'}\}=J^{\circ}_{i,j;i',j'}\ \text{and}\  i=j,\ i'=j'$;
$\triangledown$ $\colon$ $R_i R_j \in J^{\triangledown}_{i,j;i',j'}
\ \text{and}\  S_{i'} S_{j'} \notin J^{\triangledown}_{i,j;i',j'}$;
$\vartriangle$ $\colon$ $S_{i'} S_{j'} \in J^{\vartriangle}_{i,j;i',j'}
\ \text{and}\  R_i R_j \notin J^{\vartriangle}_{i,j;i',j'}$;
$\square$ $\colon$ $\{R_i R_j,S_{i'} S_{j'}\} \in J^{\square}_{i,j;i',j'}$;
}
\label{F:TS}
\end{center}
\end{figure}

Let the block $J_{i,j;i',j'}$ denote the subgraph of a joint structure $J(R,S,I)$
induced by a pair of
subsequences $\{R_i,R_{i+1},\ldots,R_j\}$ and $\{S_{i'},S_{i'+1},\ldots, S_{j'}\}$.
Given a joint structure $J(R,S,I)$, a \emph{tight} structure of $J(R,S,I)$ is
the minimal block $J_{i,j;i',j'}$ containing all the $R$-ancestors and
$S$-ancestors of any exterior arc in $J_{i,j;i',j'}$ and all the descendants
of any interior arc in $J_{i,j;i',j'}$. In the following, a tight
structure is denoted by $J^{T}_{i,j;i',j'}$. In particular, we
denote the joint structure $J(R,S,I)$ by $J^T (R,S,I)$ if $J(R,S,I)$
is a tight structure of itself. For any joint structure, there are
only four types of tight structures $J^T_{i,j;i',j'}$, that is
$\{\circ,\triangledown,\vartriangle,\square \}$, denoted by $J^{\{
\circ,\triangledown,\vartriangle,\square
 \}}_{i,j;i',j'}$, respectively, see Fig.~\ref{F:TS}.

\begin{figure}
\begin{center}
\includegraphics[width=0.65\textwidth]{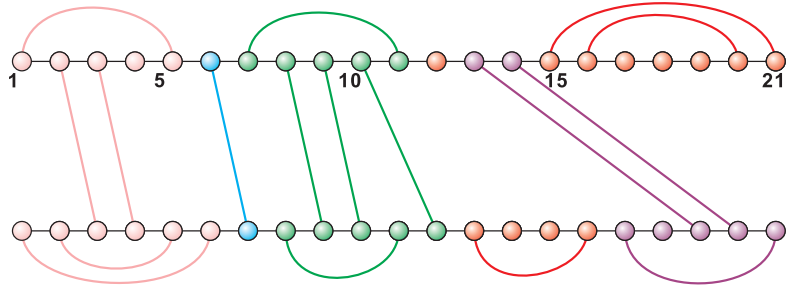}
\caption{Decomposition of joint structures. We display different secondary
segments (red) and tight structures ($\circ$ blue, $\triangledown$ green,
$\vartriangle$ purple, $\square$ pink) in which $ J_{1,21;1,21}$ decomposes.
}
\label{F:TSDecom}
\end{center}
\end{figure}

The key function of tight structures is that they are the building blocks for
the decomposition of joint structures.
\begin{proposition}\label{P:Decomposition}\cite{rip:09}
Let $J(R,S,I)$ be a joint structure. Then
\begin{enumerate}
\item any exterior arc $R_k S_{k'}$ in $J(R,S,I)$ is contained in a
      unique tight structure.
\item $J(R,S,I)$ decomposes into a unique collection of tight structures
      and maximal secondary segments.
\end{enumerate}
\end{proposition}
In Fig.~\ref{F:TSDecom} we illustrate the decomposition of a joint
structure.

\section{Refined shapes} \label{S:ReShape}

\begin{figure}[!ht]
\begin{center}
\includegraphics[width=0.8\textwidth]{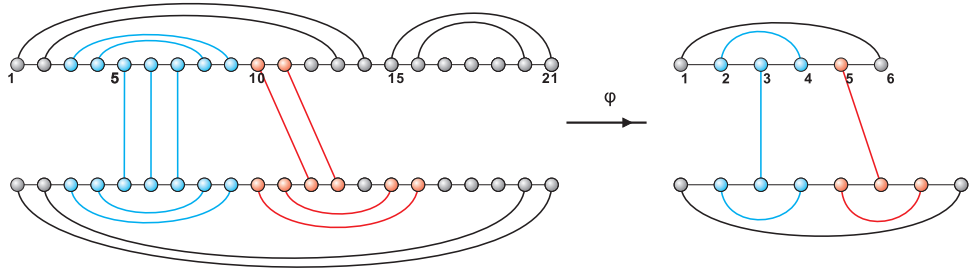}
\end{center}
\caption{\small Joint structures and their refined shapes: a $2$-canonical joint
structure with arc-length $\geq 4$ (left) is projected to its refined shape (right),
which exhibits elements of class ${\bf A}_1$ (red) and ${\bf A}_2$ (blue),
respectively.}
\label{F:RefinedShapeEx}
\end{figure}

A \emph{shape}~\cite{Li:10} is a joint structure containing no secondary segments
in which each interior stack and each exterior stack have length exactly one.
We follow the ideas in~\cite{Modular} and obtain the generating function of
joint structures via inflation of refined (colored) shapes.
Refined shapes are obtained by distinguishing two classes of exterior shape-arcs.
Each distinguished class requires its specific inflation-procedure
(see Theorem~\ref{T:RelateShape}).
Let us have a closer look at two particular classes of exterior arcs:
\begin{itemize}
\item {\bf Class $\mathbf{A}_1$:} the class of arc-pairs $(\alpha,\beta)$ where
$\alpha$ is an exterior arc with the unique interior $2$-arc $\beta$ as its
ancestor.
\item {\bf Class $\mathbf{A}_2$:} the class of arc-triples $(\alpha,\beta,\gamma)$
where $\alpha$ is an exterior arc with interior $2$-arcs $\beta$ and $\gamma$ as
ancestors.
\end{itemize}
Let $\mathcal{G}$ denote the combinatorial class of shapes. Given a joint
structure, we can obtain its shape by first removing all secondary segments
and second collapsing any stack into a single arc. That is, we have a map
$\varphi \colon {\mathcal J}\rightarrow
{\mathcal G}$. In Fig.~\ref{F:RefinedShapeEx}, we illustrate how a joint structure
is projected into its refined shape. The resulting shape exhibits elements in
class ${\bf A}_1$ as well as class ${\bf A}_2$.

Let $G(t,h,a_1,a_2)$ denote the number of shapes having total $t$ interior arcs and
$h$ exterior arcs containing $a_1$ elements of class ${\bf A}_1$ and $a_2$ elements
of class ${\bf A}_2$ and
\begin{equation*}
\mathbf{G}(x,z,u,v) = \sum_{t,h,a_1,a_2} G(t,h,a_1,a_2)\,x^t z^h u^{a_1} v^{a_2}.
\end{equation*}
We shall proceed by revisiting the notions of tight shapes, double tight shapes,
interaction segments, closed shapes and right closed shapes~\cite{rip:09}:
\begin{itemize}
\item a tight shape is tight as a structure.
    Let $\mathcal{G}^T$ denote the class of tight shapes and
    $G^T (t,h,a_1,a_2)$ denote the number of tight shapes having total
    $t$ interior arcs and $h$ exterior arcs containing $a_1$ elements of
    class ${\bf A}_1$ and $a_2$ elements of class ${\bf A}_2$,
    \begin{equation*}
    {\bf G}^T(x,z,u,v)=\sum G^T(t,h,a_1,a_2)\, x^t z^h u^{a_1} v^{a_2}.
    \end{equation*}
    Any tight shape, comes as exactly one of the four types
    $\{\circ,\triangledown,\vartriangle,\square\}$.
    The corresponding classes and generating functions are
    defined accordingly,
    $\mathcal{G}^{\{\circ,\triangledown,\vartriangle,\square\}}$
    and ${\bf
    G}^{\{\circ,\triangledown,\vartriangle,\square\}}(x,z,u,v)$
    respectively,
\item a double tight shape is a shape whose leftmost and rightmost
    blocks are tight structures. Let $\mathcal{G}^{DT}$ denote
    the class of double tight shapes and $G^{DT} (t,h,a_1,a_2)$
    denote the number of double tight shapes having total
    $t$ interior arcs and $h$ exterior arcs containing $a_1$ elements of
    class ${\bf A}_1$ and $a_2$ elements of class ${\bf A}_2$,
    \begin{equation*}
    {\bf G}^{DT}(x,z,u,v)=\sum G^{DT}(t,h,a_1,a_2)\, x^t z^h u^{a_1} v^{a_2},
    \end{equation*}
\item a closed shape is a tight shape of type $\{\triangledown,
    \vartriangle,\square\}$. Let $\mathcal{G}^C$ denote the class of
    closed shapes and $G^C (t,h,a_1,a_2)$ denote the number of closed
    shapes having total $t$ interior arcs and $h$ exterior arcs containing
    $a_1$ elements of class ${\bf A}_1$ and $a_2$ elements of class
    ${\bf A}_2$,
    \begin{equation*}
    {\bf G}^C(x,z,u,v)=\sum G^C(t,h,a_1,a_2)\, x^t z^h u^{a_1} v^{a_2},
    \end{equation*}
\item a right closed shape is a shape whose rightmost block is
    a closed shape rather than an exterior arc. Let
    $\mathcal{G}^{RC}$ denote the class of right closed shapes
    and $G^{RC} (t,h,a_1,a_2)$ denote the number of right close
    shapes having total $t$ interior arcs and $h$ exterior arcs containing
    $a_1$ elements of class ${\bf A}_1$ and $a_2$ elements of class ${\bf A}_2$,
    \begin{equation*}
    {\bf G}^{RC}(x,z,u,v)=\sum G^{RC}(t,h,a_1,a_2)\, x^t z^h u^{a_1} v^{a_2},
    \end{equation*}
\item in a shape, an interaction segment is an empty structure or a tight
    structure of type $\circ$ (an exterior arc). We denote the
    class of interaction segment by $\mathcal{I}$ and the associated
    generating function by ${\bf I}(x,z,u,v)$. Obviously, ${\bf
    I}(x,z,u,v)=1+z$.
\end{itemize}

\begin{lemma}\label{L:Shape}
The generating function ${\bf G}(x,z,u,v)$ of refined shapes satisfies
\begin{equation}\label{E:Shape}
{\bf A}(x,z){\bf G}(x,z,u,v)^2+ {\bf B}(x,z,u,v){\bf G}(x,z,u,v)+
{\bf C}(x,z)=0,
\end{equation}
where
\begin{equation}\label{E:erni}
\begin{aligned}
{\bf A}(x,z)  &= x(x+2)(z+1),\\
{\bf B}(x,z,u,v)  &= -\left(x(x+2)(z+1)^2 + (x+1)^2 -(2 \,x \, u + x^2\, v )
z(z+1) \right),\\
{\bf C}(x,z)  &= (x+1)^2 (z+1) .
\end{aligned}
\end{equation}
Explicitly,
\begin{equation}\label{E:ShapeExp}
{\bf G}(x,z,u,v) = \frac{-\mathbf{B}(x,z,u,v)-\sqrt{\mathbf{B}(x,z,u,v)^2-4 \,
\mathbf{A}(x,z) \mathbf{C}(x,z)}}{2\,\mathbf{A}(x,z)}.
\end{equation}
\end{lemma}
\begin{figure}
\begin{center}
\includegraphics[width=0.45\textwidth]{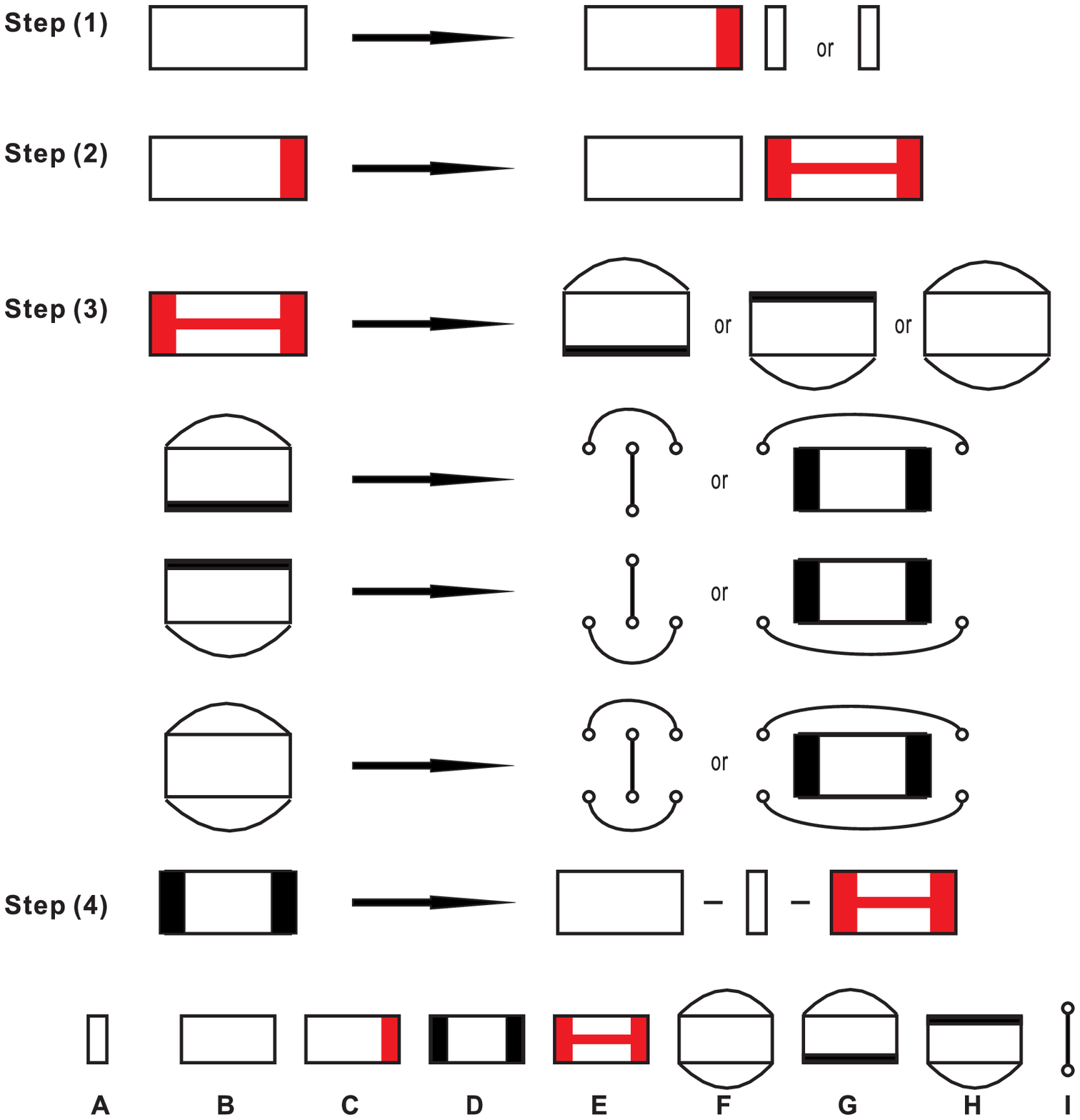}
\caption{The shape-grammar. The notations of structural components are
explained in the panel below.
{\bf{A}}: interaction segment;
{\bf{B}}: arbitrary shape $G(R,S,I)$; {\bf{C}}: right closed shape
$G^{RC}(R,S,I)$; {\bf{D}}: double tight shape $G^{DT}(R,S,I)$;
{\bf{E}}: closed shape $G^{C}(R,S,I)$; {\bf{F}}: type $\square$
tight shape $G^{\square}(R,S,I)$; {\bf{G}}: type $\triangledown
$ tight shape $G^{\triangledown }(R,S,I)$; {\bf{H}}: type
$\vartriangle$ tight shape $G^{\vartriangle}(R,S,I)$; {\bf{I}}:
type $\circ$ tight shape $G^{\circ}(R,S,I)$; }
\label{F:ShapeDecomGram}
\end{center}
\end{figure}
\begin{proof}
Proposition \ref{P:Decomposition} implies that any shape can be
decomposed into a unique collection of tight shapes. Furthermore,
each shape can be decomposed into a unique collection of closed shapes
and exterior arcs. We decompose a shape in four steps, see
Fig.~\ref{F:ShapeDecomGram}. We translate each decomposition step into
the construction of combinatorial classes.\\
{\bf Step (1):} decomposition into a right closed shape and
its rightmost interaction segment, i.e.
$$
\mathcal{G}=\mathcal{G}^{RC}\times\mathcal{I}+ \mathcal{I}.
$$
and Proposition~\ref{P:Symbolic} implies
\begin{equation}\label{E:Shape-1}
{\bf G}(x,z,u,v)={\bf G}^{RC}(x,z,u,v)\cdot {\bf I}(x,z,u,v)+{\bf I}(x,z,u,v).
\end{equation}
{\bf Step (2):} splitting of the rightmost closed shape in a right closed shape
\[
\mathcal{G}^{RC}=\mathcal{G}\times \mathcal{G}^{C},
\]
whence
\begin{equation}\label{E:Shape-2}
{\bf G}^{RC}(x,z,u,v)={\bf G}(x,z,u,v)\cdot {\bf G}^{C}(x,z,u,v).
\end{equation}
{\bf Step (3):} type-depended decomposition of a closed shape.
\begin{eqnarray*}
\mathcal{G}^{C}               &=&
\mathcal{G}^{\triangledown}+\mathcal{G}^{\vartriangle}
+\mathcal{G}^{\square}\\
\mathcal{G}^{\triangledown}   &=&
(\mathcal{Z},\mathcal{Z},\mathcal{Z},\mathcal{E}) +
(\mathcal{Z},\mathcal{E},\mathcal{E},\mathcal{E}) \times \mathcal{G}^{DT}\\
\mathcal{G}^{\vartriangle}   &=&
(\mathcal{Z},\mathcal{Z},\mathcal{Z},\mathcal{E})+
(\mathcal{Z},\mathcal{E},\mathcal{E},\mathcal{E}) \times \mathcal{G}^{DT}\\
\mathcal{G}^{\square}       &=&
(\mathcal{Z}\times \mathcal{Z},\mathcal{Z},\mathcal{E},\mathcal{Z})+
(\mathcal{Z}\times \mathcal{Z},\mathcal{E},\mathcal{E}) \times \mathcal{G}^{DT}.
\end{eqnarray*}
We therefore have
\begin{equation}\label{E:Shape-3}
\begin{aligned}
{\bf G}^{C}(x,z,u,v)                  &=  {\bf G}^{\triangledown}(x,z,u,v)+
{\bf G}^{\vartriangle}(x,z,u,v)+{\bf G}^{\square}(x,z,u,v)\\
{\bf G}^{\triangledown}(x,z,u,v)      &=  x\, z \, u +  x \, {\bf G}^{DT}(x,z,u,v)\\
{\bf G}^{\vartriangle}(x,z,u,v)       &=  x \, z\, u +  x \, {\bf G}^{DT}(x,z,u,v)\\
{\bf G}^{\square}(x,z,u,v)            &=  x^2 \,z\,v +  x^2 \, {\bf G}^{DT}(x,z,u,v).
\end{aligned}
\end{equation}
{\bf Step (4):} obtaining double tight shapes of Step (3) by excluding
the class of interaction segment and the class of closed shapes, i.e.
\[
\mathcal{G}^{DT}=\mathcal{G}-\mathcal{I}-\mathcal{G}^C,
\]
whence
\begin{equation}\label{E:Shape-4}
{\bf G}^{DT}(x,z,u,v)={\bf G}(x,z,u,v) - {\bf I}(x,z,u,v) -{\bf G}^{C}(x,z,u,v).
\end{equation}
Solving equations~(\ref{E:Shape-1})--(\ref{E:Shape-4}), leads to (1) the functional
equation~(\ref{E:Shape}) and (2) eq.~(\ref{E:erni}).
This quadratic equation together with the initial condition
${\bf G}(0,0,0,0)=1$, implies eq.~(\ref{E:ShapeExp}).
\end{proof}


\section{The generating function of joint structures}\label{S:inflate}

Let $\mathcal{H}_{\sigma}$ denote the class of $\sigma$-canonical joint structures
with arc-length $\geq \sigma+2$. Let $H_{\sigma}(s)$ denote the number of joint
structures in $\mathcal{H}_{\sigma}$ with total $s$ vertices having the generating
function
\begin{equation*}
\mathbf{H}_{\sigma}(x) = \sum H_{\sigma}(s) \, x^s.
\end{equation*}
We are now in position to compute the generating function $\mathbf{H}_{\sigma}(x)$.
Our strategy is inflating the shapes with specific inflation on specific exterior
arcs.
\begin{theorem}\label{T:RelateShape}
For any $\sigma \geq 1$, $\mathbf{H}_{\sigma}(x)$ is a power series and
\begin{equation}\label{E:RelateShape}
\mathbf{H}_{\sigma}(x) =
{\bf T}_{\sigma}^{[\sigma+2]}(x)^2 \,
{\bf G}\!\left(\eta,\eta,\eta_1 ,\eta_2 \right),
\end{equation}
where
\begin{eqnarray*}
\eta     &=&  \frac{x^{2\sigma}\,
{\bf T}_{\sigma}^{[\sigma+2]}(x)^2}{1-x^2-x^{2\sigma}
({\bf T}_{\sigma}^{[\sigma+2]}(x)^2-1)},\\
\eta_1   &=&  \frac{-1+ x^2 -x^{2\sigma}+(1+x^{2\sigma})
           {\bf T}_{\sigma}^{[\sigma+2]}(x)^2}{{\bf T}_{\sigma}^{[\sigma+2]}(x)^2},\\
\eta_2   &=&  \frac{1-x^2+x^{2\sigma}+(-2+2x^2-3x^{2\sigma})
              {\bf T}_{\sigma}^{[\sigma+2]}(x)^2+(1+2x^{2\sigma})
              {\bf T}_{\sigma}^{[\sigma+2]}(x)^4}{{\bf T}_{\sigma}^{[\sigma+2]}(x)^4}.
\end{eqnarray*}
\end{theorem}

\begin{proof}
Let ${\mathcal G}(t,h,a_1,a_2)$ denote the class of shapes having total $t$
interior arcs and $h$ exterior arcs containing $a_1$ elements of class ${\bf A}_1$
and $a_2$ elements of class ${\bf A}_2$. For any joint structure in
$\mathcal{H}_{\sigma}$, we obtain a shape in $\mathcal{G}$ as follows:
\begin{enumerate}
\item remove all secondary segments,
\item collapse each interior stack into one interior arc and each exterior
stack into one exterior arc.
\end{enumerate}
Then we have the surjective map
\begin{equation*}
\varphi\colon \mathcal{H}_{\sigma} \rightarrow {\mathcal G}.
\end{equation*}
Indeed, for any shape $\gamma$ in $\mathcal{G}$, we can construct
$\sigma$-canonical joint structures with arc-length $\geq \sigma+2$.
$\varphi\colon \mathcal{H}_{\sigma} \rightarrow {\mathcal G}$, induces the
partition $\mathcal{H}_{\sigma}=\dot\cup_\gamma\varphi^{-1}(\gamma)$,
whence
\begin{equation}
\mathbf{H}_{\sigma}(x) = \sum_{\gamma\in\,{\mathcal G}}
\mathbf{H}_\gamma(x),
\end{equation}
where ${\bf H}_\gamma(x)$ denotes the generating function of joint structures
having shape $\gamma$.
We proceed by computing the generating function
$\mathbf{H}_\gamma(x)$. We will construct
$\mathbf{H}_\gamma(x)$ via simpler combinatorial classes as
building blocks considering stems, stacks, induced stacks, interior arcs,
exterior arcs and
secondary segments. We inflate a shape $\gamma$ in ${\mathcal G}(t,h,a_1,a_2)$
to a joint structure in $\mathcal{H}_{\sigma}$ in four steps.\\
{\bf Step I:} we inflate any interior arc in $\gamma$ to a stack of
size at least $\sigma$ and subsequently add additional stacks. The
latter are called induced stacks and have to be separated by means
of inserting secondary segments, see Fig.~\ref{F:RelateJG-1}.
\begin{figure}
\includegraphics[width=0.9\textwidth]{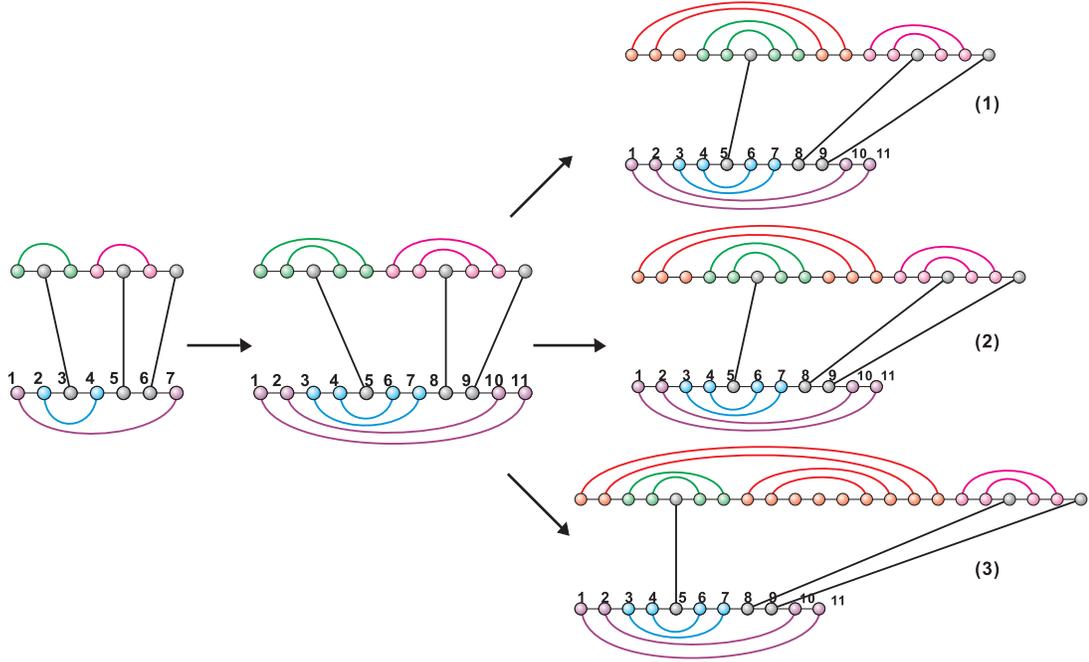}
\caption{Step I: a shape (left) is inflated to a joint structure with arc-length
$\geq 2$ and interior stack-length $\geq 2$. Each interior arc in the shape
is first inflated to a stack of size at least two (middle) and second inflated by
adding one induced stack of size two (right). Note that there are three ways
to insert the secondary segments to separate the induced stacks (red).}
\label{F:RelateJG-1}
\end{figure}
Note that during this first inflation step no secondary segments,
other than those necessary for separating the nested stacks are
inserted. We generate
\begin{itemize}
\item secondary segments $\mathcal{T}_{\sigma}^{[\sigma+2]}$ having arc-length
$\geq \sigma+2$ and stack-length $\geq \sigma$ with generating function
${\bf T}_{\sigma}^{[\sigma+2]}(x)$,
\item interior arcs $\mathcal{R}$ with generating function ${\bf R}(x)=x^2$,
\item stacks, i.e.~pairs consisting of the minimal sequence of arcs
$\mathcal{R}^\sigma$ and an arbitrary extension consisting of
arcs of arbitrary finite length
\begin{equation*}
\mathcal{K}_{\sigma}=\mathcal{R}^{\sigma}\times\textsc{Seq}\left(\mathcal{R}\right),
\end{equation*}
having the generating function
\begin{eqnarray*}
\mathbf{K}_\sigma(x) & = & x^{2\sigma}\cdot \frac{1}{1-x^2},
\end{eqnarray*}
\item induced stacks, i.e.~stacks together with at least one secondary segment on
either or both of its sides,
\begin{equation*}
\mathcal{N}_{\sigma}=\mathcal{K}_{\sigma}\times\left(
(\mathcal{T}_{\sigma}^{[\sigma+2]})^2 -1 \right),
\end{equation*}
having the generating function
\begin{equation*}
\mathbf{N}_\sigma(x)=\frac{x^{2\sigma}}{1-x^2}
\left( {\bf T}_{\sigma}^{[\sigma+2]}(x)^2-1\right),
\end{equation*}
\item stems, that is pairs consisting of stacks $\mathcal{K}_\sigma$
and an arbitrarily long sequence of induced stacks
\begin{equation*}
\mathcal{M}_\sigma=\mathcal{K}_{\sigma} \times
\textsc{Seq}\left(\mathcal{N}_{\sigma}\right),
\end{equation*}
having the generating function
\begin{eqnarray*}
\mathbf{M}_\sigma(x)=\frac{\mathbf{K}_\sigma(x)}{1-\mathbf{N}_\sigma(x)}=
\frac{\frac{x^{2\sigma}}{1-x^2}}
{1-\frac{x^{2\sigma}}{1-x^2}\left( {\bf T}_{\sigma}^{[\sigma+2]}(x)^2-1 \right)}.
\end{eqnarray*}
\end{itemize}
Note that we inflate both, top and bottom sequences. The
corresponding generating function is ${\bf M}_{\sigma}(x)^{t}$.

{\bf Step II:} we inflate any exterior arc in $\gamma$, but not as an element of
classes ${\bf A}_1$ or ${\bf A}_2$, to an exterior stack of size at least $\sigma$
and subsequently add additional exterior stacks. The latter are called induced exterior
stacks and have to be separated by means of inserting secondary
segments, see Fig.~\ref{F:RelateJG-2}.
\begin{figure}
\includegraphics[width=0.9\textwidth]{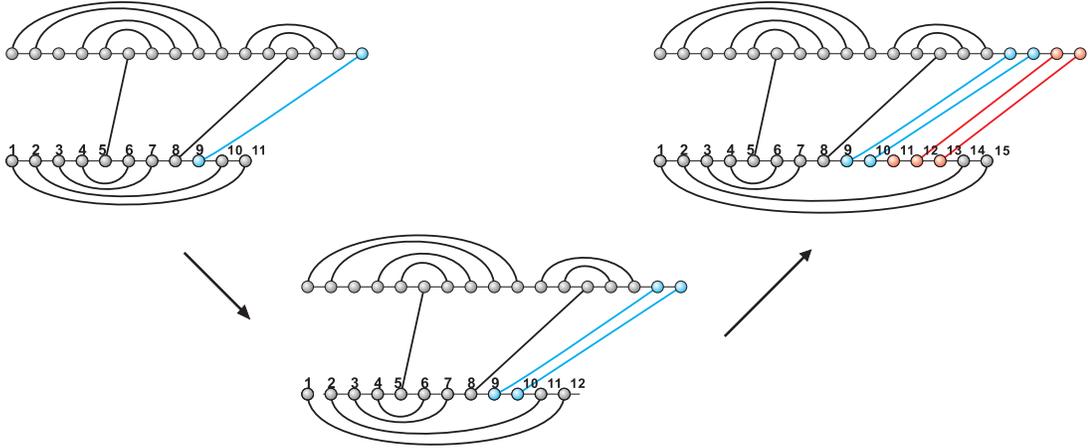}
\caption{Step II: a joint structure (top-left) obtained in {\sf (1)} in
Fig.~\ref{F:RelateJG-1} is inflated to a new joint structure.
Each exterior arc, not contained in classes
${\bf A}_1$ or ${\bf A}_2$, is first inflated to
an exterior stack of size at least two (bottom) and second inflated by adding
one exterior induced stack of size two (top-right). Note that just one of the three
possible ways of inserting the secondary segments in order to separate the induced
exterior stacks (red) is displayed.
} \label{F:RelateJG-2}
\end{figure}
Note that during this exterior-arc inflation no secondary
segments, other than those necessary for separating the stacks are
inserted. We generate
\begin{itemize}
\item exterior arc $\mathcal{R}_0$ having the generating function ${\bf R}_0= x^2$,
\item exterior stacks, i.e.~pairs consisting of the minimal sequence of exterior arcs
$\mathcal{R}_0^\sigma$ and an arbitrary extension consisting of
exterior arcs of arbitrary finite length
\begin{equation*}
\mathcal{K}^{\ast}_{\sigma}=
\mathcal{R}_0^{\sigma}\times\textsc{Seq}\left(\mathcal{R}_0\right),
\end{equation*}
having the generating function
\begin{eqnarray*}
\mathbf{K}^{\ast}_{\sigma}(x) & = & x^{2 \sigma}\cdot \frac{1}{1-x^2},
\end{eqnarray*}
\item induced exterior stacks, i.e.~stacks together with at least one secondary
segment on either or both its sides,
\begin{equation*}
\mathcal{N}^{\ast}_{\sigma}=\mathcal{K}^{\ast}_{\sigma}\times
\left( (\mathcal{T}_{\sigma}^{[\sigma+2]})^2 -1 \right),
\end{equation*}
having the generating function
\begin{equation*}
\mathbf{N}^{\ast}_\sigma(x) =\frac{x^{2 \sigma}}{1-x^2 }
\left( {\bf T}_{\sigma}^{[\sigma+2]}(x)^2-1\right),
\end{equation*}
\item exterior stems, that is pairs consisting of exterior stacks
$\mathcal{K}^{\ast}_\sigma$ and an arbitrarily long sequence of induced exterior
stacks
\begin{equation*}
\mathcal{M}^{\ast}_\sigma=\mathcal{K}^{\ast}_{\sigma} \times
                                \textsc{Seq}\left(\mathcal{N}^{\ast}_{\sigma}\right),
\end{equation*}
having the generating function
\begin{eqnarray*}
\mathbf{M}^{\ast}_\sigma(x)=\frac{\mathbf{K}^{\ast}_\sigma(x)}
{1-\mathbf{N}^{\ast}_\sigma(x)}=
\frac{\frac{x^{2 \sigma}}{1-x^2 }}
{1-\frac{x^{2 \sigma}}{1-x^2 } \left( {\bf T}_{\sigma}^{[\sigma+2]}(x)^2-1\right)}.
\end{eqnarray*}
\end{itemize}
Inflating all exterior arcs that are not contained in classes
${\bf A}_1$ or ${\bf A}_2$, we obtain $({\bf M}^{\ast}_{\sigma}(x))^{h-a_1-a_2}$.

{\bf Step III:} We inflate exterior arcs contained in classes ${\bf A}_1$ and
${\bf A}_2$ by inserting additional secondary segments at positions between
the exterior arc and interior $2$-arc, see Fig.~\ref{F:RelateJG-3}.
In contrast to Step II, specific ``unwanted'' scenarios are excluded.
\begin{figure}
\includegraphics[width=0.9\textwidth]{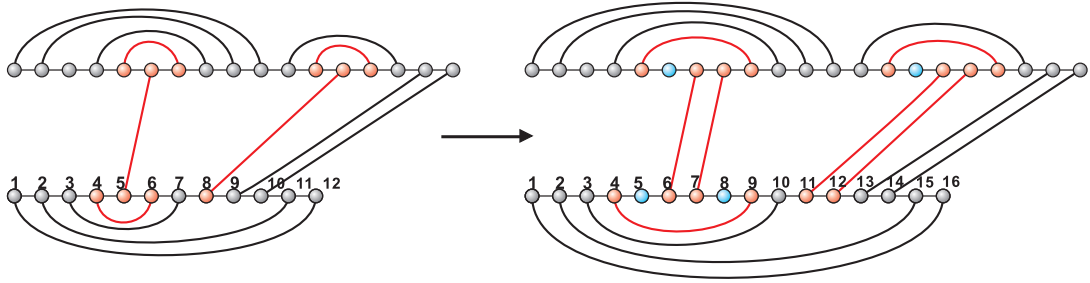}
\caption{Step III: a joint structure (left) obtained in the bottom in
Fig.~\ref{F:RelateJG-2} is inflated to a new joint structure with arc-length
$\geq 4$ (right). Each exterior arc in classes ${\bf A}_1$ or ${\bf A}_2$ is
inflated to an exterior stack of size at least two (red) and additional
secondary segments (blue) are inserted at the positions between the exterior
arc and interior $2$-arc.} \label{F:RelateJG-3}
\end{figure}
We generate
\begin{itemize}
\item {\bf Class $\mathbf{A}_1$:}
Excluding the case where the exterior arc is inflated to an exterior stack of length
$\sigma$ and no additional secondary segment is inserted at the position between the
exterior arc and interior $2$-arc, see Fig.~\ref{F:Class1}, we arrive at
    \begin{equation*}
    \mathcal{M}^{\ast}_\sigma \times
    ({\mathcal{T}}_{\sigma}^{[\sigma+2]})^2-\mathcal{R}_0^{\sigma},
    \end{equation*}
    having the generating function
    \begin{equation*}
    \mathbf{M}^{\ast}_\sigma(x) {\mathbf{T}}_{\sigma}^{[\sigma+2]}(x)^2 - x^{2 \sigma}.
    \end{equation*}
\begin{figure}
\begin{center}
\includegraphics[width=0.4\textwidth]{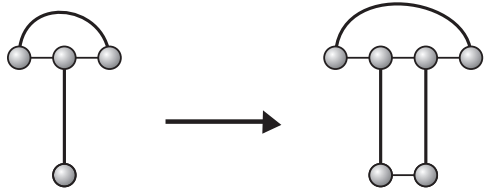}
\caption{{\bf ``Bad'' scenario for class ${\bf A}_1$:} In an element of class
$\mathbf{A}_1$ (left), the exterior arc is inflated to an exterior stack of length
$2$ and no additional secondary segment is inserted at the position between the
exterior arc and interior $2$-arc, leading to an interior arc having arc-length
$<4$ (right).} \label{F:Class1}
\end{center}
\end{figure}
\item {\bf Class $\mathbf{A}_2$:} There are three scenarios which create an
interior arc of arc-length $<\sigma +2$, see Fig.~\ref{F:Class2}:
    \begin{itemize}
    \item  the exterior arc is inflated to an exterior stack of length $\sigma$
and no additional secondary segment is inserted at the positions in both top and
bottom sequences, resulting in both interior arcs having arc-length $<\sigma +2$,
    \item  the exterior arc is inflated to an exterior stack of length $\sigma$
and additional secondary segment is inserted at the positions only in the bottom
sequence, resulting in an interior arc in the top having arc-length $<\sigma +2$,
    \item  the exterior arc is inflated to an exterior stack of length $\sigma$
and additional secondary segment is inserted at the positions only in the top
sequence, resulting in an interior arc in the bottom having arc-length $<\sigma +2$.
    \end{itemize}
\begin{figure}
\begin{center}
\includegraphics[width=0.5\textwidth]{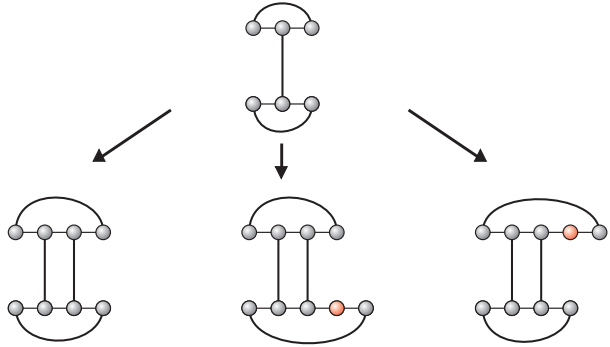}
\caption{{\bf ``Bad'' scenarios for class ${\bf A}_2$:} In an element of class
$\mathbf{A}_2$ (top), there are three scenarios leading to interior arcs in both
top and bottom having arc-length $<4$ (bottom-left), the interior arc only in the
 top having arc-length $<4$ (bottom-middle), the interior arc only in the bottom
having arc-length $<4$ (bottom-right).
} \label{F:Class2}
\end{center}
\end{figure}
Excluding these three scenarios, we obtain
\begin{equation*}
\mathcal{M}^{\ast}_\sigma \times ({\mathcal{T}}_{\sigma}^{[\sigma+2]})^4 -2\,
\mathcal{R}_0^{\sigma}\times (({\mathcal{T}}_{\sigma}^{[\sigma+2]})^2-
\mathcal{E} )-\mathcal{R}_0^{\sigma},
\end{equation*}
having the generating function
\begin{eqnarray*}
&&  \mathbf{M}^{\ast}_\sigma(x) {\mathbf{T}}_{\sigma}^{[\sigma+2]}(x)^4-
x^{2 \sigma}({\bf T}_{\sigma}^{[\sigma+2]}(x)^2-1) -
x^{2 \sigma}({\bf T}_{\sigma}^{[\sigma+2]}(x)^2-1) - x^{2 \sigma}\\
&=& \mathbf{M}^{\ast}_\sigma(x) {\bf T}_{\sigma}^{[\sigma+2]}(x)^4 -
x^{2 \sigma}( 2{\bf T}_{\sigma}^{[\sigma+2]}(x)^2-1).
\end{eqnarray*}
\end{itemize}
Applying Step III for each ${\bf A}_1$- and ${\bf A}_2$-element, we derive
\begin{eqnarray*}
\left(\mathbf{M}^{\ast}_\sigma(x) {\bf T}_{\sigma}^{[\sigma+2]}(x)^2 -
x^{2 \sigma}\right)^{a_1}
\left(\mathbf{M}^{\ast}_\sigma(x) {\bf T}_{\sigma}^{[\sigma+2]}(x)^4 -
x^{2 \sigma}( 2{\bf T}_{\sigma}^{[\sigma+2]}(x)^2-1) \right)^{a_2}.
\end{eqnarray*}

{\bf Step IV:} Here we insert additional secondary segments at the
remaining  $(2t+2h+2-2 a_1-4a_2)$ positions, see Fig.~\ref{F:RelateJG-4}.
Formally, this fourth inflation is expressed via the combinatorial class
\begin{equation*}
(\mathcal{T}_{\sigma}^{[\sigma+2]})^{2t+2h+2-2 a_1-4a_2} ,
\end{equation*}
with the generating function $({\bf T}_{\sigma}^{[\sigma+2]}(x))^{2t+2h+2-2 a_1-4a_2}$.
\begin{figure}
\includegraphics[width=0.9\textwidth]{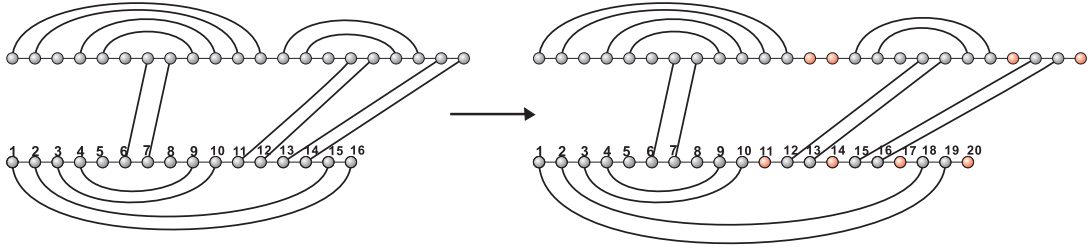}
\caption{Step IV: a joint structure (left) obtained in
Fig.~\ref{F:RelateJG-3} is inflated to a new joint structure in ${\mathcal H}_{2}$
(right) by inserting secondary segments (red).} \label{F:RelateJG-4}
\end{figure}

Combining Steps I -- IV, we arrive at
\begin{eqnarray*}
{\bf H}_\gamma(x)&=&
{\bf M}_{\sigma}(x)^{t}\,({\bf M}^{\ast}_{\sigma}(x))^{h-a_1-a_2}
\left(\mathbf{M}^{\ast}_\sigma(x) {\bf T}_{\sigma}^{[\sigma+2]}(x)^2 - x^{2 \sigma}
\right)^{a_1}\\
&&\times\left(\mathbf{M}^{\ast}_\sigma(x) {\bf T}_{\sigma}^{[\sigma+2]}(x)^4 -
x^{2 \sigma}( 2{\bf T}_{\sigma}^{[\sigma+2]}(x)^2-1) \right)^{a_2}
({\bf T}_{\sigma}^{[\sigma+2]}(x))^{2t+2h+2-2 a_1-4a_2}  \\
&=&({\bf T}_{\sigma}^{[\sigma+2]}(x))^{2} (\eta)^{t} (\eta)^h (\eta_1)^{a_1}
(\eta_2)^{a_2},
\end{eqnarray*}
where
\begin{eqnarray*}
\eta     &=&  \frac{x^{2\sigma}\,
{\bf T}_{\sigma}^{[\sigma+2]}(x)^2}{1-x^2-x^{2\sigma}
({\bf T}_{\sigma}^{[\sigma+2]}(x)^2-1)},\\
\eta_1      &=&  \frac{\eta-x^{2\sigma}}{\eta}\\
            &=&  \frac{-1+ x^2 -x^{2\sigma}+(1+x^{2\sigma}){\bf T}_{\sigma}^{[\sigma+2]}(x)^2}{{\bf T}_{\sigma}^{[\sigma+2]}(x)^2},\\
\eta_2      &=&  \frac{\eta {\bf T}_{\sigma}^{[\sigma+2]}(x)^2 -x^{2\sigma}(2{\bf T}_{\sigma}^{[\sigma+2]}(x)^2-1)}{\eta {\bf T}_{\sigma}^{[\sigma+2]}(x)^2}\\
            &=&  \frac{1-x^2+x^{2\sigma}+(-2+2x^2-3x^{2\sigma}){\bf T}_{\sigma}^{[\sigma+2]}(x)^2+(1+2x^{2\sigma})
            {\bf T}_{\sigma}^{[\sigma+2]}(x)^4}{{\bf T}_{\sigma}^{[\sigma+2]}(x)^4}.
\end{eqnarray*}
In view of ${\bf T}_{\sigma}^{[\sigma+2]}(0)\neq 0$ the inverse
$[{\bf T}_{\sigma}^{[\sigma+2]}(x)]^{-1}$ exists and accordingly $\eta$, $\eta_1$ and
$\eta_2$ are welldefined.
Since for any $\gamma,\gamma_1\in {\mathcal G}(t,h,a_1,a_2)$ we have
$\mathbf{H}_\gamma(x)=\mathbf{H}_{\gamma_1}(x)$, we derive
\begin{equation*}
{\bf H}_{\sigma}(x) = \sum_{\gamma\in\,{\mathcal G}}
\mathbf{H}_\gamma(x) =
\sum_{(t,h,a_1,a_2) \atop \gamma\in\,{\mathcal G}(t,h,a_1,a_2)}
G(t,h,a_1,a_2)\,\mathbf{H}_\gamma(x).
\end{equation*}
Using  ${\bf G}(x,z,u,v)= \sum G(t,h,a_1,a_2)\, x^t z^h u^{a_1} v^{a_2}$,
we arrive at
\begin{equation*}
\mathbf{H}_{\sigma}(x) =
{\bf T}_{\sigma}^{[\sigma+2]}(x)^2 \,
{\bf G}\!\left(\eta,\eta,\eta_1 ,\eta_2 \right).
\end{equation*}
It remains to verify that $H_\sigma(x)$ is indeed a power series, which follows from
the fact that the constant coefficients of $\eta$, $\eta_1$ and $\eta_2$, regarded
as formal power series, are zero.
\end{proof}


\section{Asymptotic Enumeration}\label{S:asym}

In this section, we derive simple formulas for the number of joint structures
in the limit of long sequences.
\begin{theorem}\label{T:AsymtoticJS}
For $\sigma \geq 1$, $\mathbf{H}_{\sigma}(x)$ is algebraic and we have
\begin{equation}\label{E:JointAsymptotic}
H_{\sigma}(s) \, \sim \, c_{\sigma} \, s^{-\frac{3}{2}} \,
\left(\kappa_{\sigma}^{-1}\right)^s ,\quad \text{for some
$c_{\sigma}$,}\qquad
\end{equation}
where $\kappa_{\sigma}$ is the minimal, positive real solution of
the equation $\mathbf{Q}(x,{\bf T}_\sigma^{[\sigma+2]}(x))=0$,
see Table~\ref{Tab:asymJS}, where
\begin{eqnarray*}
\mathbf{Q}(x,y)   &=&
\left(1-x^2+x^{2 \sigma }\right)^4+2 x^{2 \sigma } \left(1-x^2+x^{2 \sigma }\right)^2
                    \left(-3+3 x^2-x^{2 \sigma } +x^{4 \sigma }-2 x^{2+2 \sigma } \right)y^2\\
                && x^{4 \sigma } \left(3-6 x^2+3 x^4+10 x^{2 \sigma }+13 x^{4 \sigma }+6 x^{6 \sigma }+x^{8 \sigma } \right. \\
                && \left. -14 x^{2+2 \sigma }+4 x^{4+2 \sigma }-14 x^{2+4 \sigma }+4 x^{4+4 \sigma }-4 x^{2+6 \sigma }\right)
                y^4\\
                &&-2 x^{6 \sigma } \left(1-x^2+3 x^{2 \sigma }+x^{4 \sigma }-2 x^{2+2 \sigma }\right)
                y^6+x^{8 \sigma } y^8.
\end{eqnarray*}
In particular, we have $c_{1}\approx 1.38629$ and $c_{2}\approx 3.51610$.
\end{theorem}
\begin{proof}
Set $\mathbf{G}_1(x)={\bf G}\!\left(\eta,\eta,\eta_1 ,\eta_2
\right)$ and $L= \mathbb{C}(x)[ {\bf T}_{\sigma}^{[\sigma+2]}(x)]$.
Combining eq.~(\ref{E:Shape}) in Lemma~\ref{L:Shape} and
eq.~(\ref{E:RelateShape}) in Theorem~\ref{T:RelateShape}, we compute
\begin{eqnarray*}
{\bf H}_{\sigma}(x) &=& {\bf T}_{\sigma}^{[\sigma+2]}(x)^2 \,
{\bf G}\!\left(\eta,\eta,\eta_1 ,\eta_2 \right)\\
                    &=& {\bf T}_{\sigma}^{[\sigma+2]}(x)^2 \, \mathbf{G}_1(x),
\end{eqnarray*}
where $\mathbf{G}_1(x)$ satisfies the quadratic equation
\begin{equation}\label{E:Quadratic}
{\bf A}(\eta,\eta){\bf G}_1(x)^2+ {\bf B}(\eta,\eta,\eta_1,\eta_2){\bf G}_1(x)+
{\bf C}(\eta,\eta)=0.
\end{equation}
Note that ${\bf A}(\eta,\eta)$, $ {\bf B}(\eta,\eta,\eta_1,\eta_2)$
and ${\bf C}(\eta,\eta)$ are elements of the quadratic field
extension $L/\mathbb{C}(x)$. Thus eq.~(\ref{E:Quadratic}) implies
that $\mathbf{G}_1(x)$ is algebraic over $L$, that is, the field
extension $ L[\mathbf{G}_1(x)]/L$ is finite, whence the extension $
L[\mathbf{G}_1(x)]/\mathbb{C}(x)$ is finite. Therefore
$\mathbf{G}_1(x)$ is algebraic over $\mathbb{C}(x)$. Clearly this
implies that $\mathbf{H}_{\sigma}(x)$ is algebraic over
$\mathbb{C}(x)$ and in particular $D$-finite~\cite{Stanley}.
Pringsheim's Theorem~\cite{Tichmarsh:39} guarantees that
$\mathbf{H}_{\sigma}(x)$ has a dominant real positive singularity
$\kappa_{\sigma}$. We verify by explicit computation that for $1\leq
\sigma \leq 5$, the singularity $\kappa_{\sigma}$ is the unique,
minimal, positive real solution of the equation
$\mathbf{Q}(x,{\bf T}_\sigma^{[\sigma+2]}(x))=0$
and a branch-point singularity of the square root. We list the
values of $\kappa_{\sigma}^{-1}$ in Table~\ref{Tab:asymJS}.
Accordingly, at $\kappa_{\sigma}$, $\mathbf{H}_{\sigma}(x)$
coincides with its singular expansion and is given by
\begin{equation*}
\mathbf{H}_{\sigma}(x)=
h_0+h_1(\kappa_{\sigma} -x )^{\frac{1}{2}}+ O((\kappa_{\sigma}-x)).
\end{equation*}
Using Theorem~\ref{T:transfer1} and Theorem~\ref{T:transfer2}, we
arrive at
\begin{equation*}
H_{\sigma}(s) \, \sim \, \frac{h_1
(\kappa_{\sigma})^{\frac{1}{2}}}{\Gamma (-\frac{1}{2})} \,
s^{-\frac{3}{2}} \, \left(\kappa_{\sigma}\right)^{-s}.
\end{equation*}
Setting $c_{\sigma} = \frac{h_1
(\kappa_{\sigma})^{\frac{1}{2}}}{\Gamma (-\frac{1}{2})}$, we compute
$c_{1}\approx 1.38629$ and $c_{2}\approx 3.51610$, completing the
proof.
\end{proof}

\begin{figure*}
\begin{center}
\includegraphics[width=0.5\textwidth]{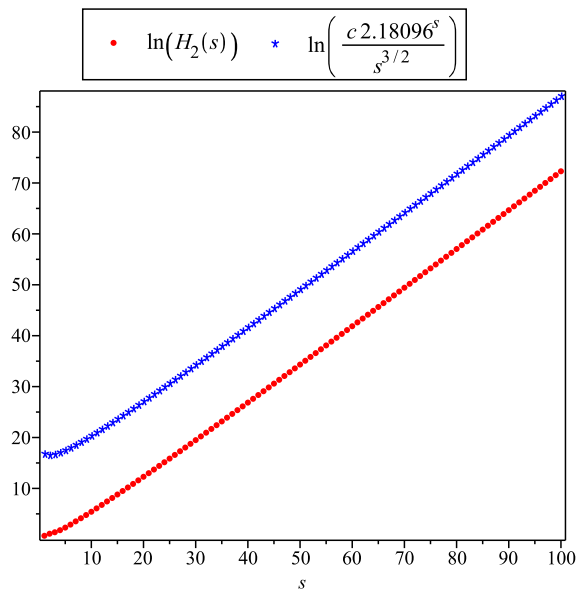}
\end{center}
\caption{Exact enumeration versus asymptotic formula.
We contrast the numbers of $2$-canonical joint structures with arc-length
$\geq 4$ ($H_2(s)$) versus $c \, s^{-\frac{3}{2}} \, 2.18096^s$.
For representational purposes we separate the curves via setting the
constant $c=10^7$.}\label{F:AsymJoint}
\end{figure*}

\begin{table}
\begin{center}
\begin{tabular}{cccccc}
\hline\noalign{\smallskip}
\small$\sigma$                &\small$1$   & \small$2$ & \small$3$  & \small$4$ &
\small $5$    \\
\noalign{\smallskip}\hline\noalign{\smallskip}
\small$\kappa_{\sigma}^{-1}$  &\small$3.30027$ & \small$2.18096$ & \small
$1.82912$ &\small $1.65183$ & \small $1.54322$ \\
\noalign{\smallskip}\hline
\end{tabular}
\end{center}
\caption{Exponential growth rates $\kappa_{\sigma}^{-1}$ for $\sigma$-canonical
joint structures with arc-length $\geq \sigma+2$.}
\label{Tab:asymJS}
\end{table}
In Fig.~\ref{F:AsymJoint}, we showcase the quality of the asymptotic formula
for $\sigma=2$ and arc-length four, implied by Theorem~\ref{T:AsymtoticJS}.

\section{Discussion}\label{S:discussion}


In this paper we analyzed the biologically relevant class of
canonical joint structures having arc-length greater than or equal
to four. While it is straightforward to derive the generating
function of joint structures from the (eleven) recursion relations
of the original {\tt rip}-grammar (implied by
Proposition~\ref{P:Decomposition})~\cite{rip:09} the generating
function obtained this way would be ``impossible'' to write down.
This approach would be neither suitable for deriving any asymptotic
formulas nor would it allows us to deal with specific stack-length
conditions. Therefore we do not use the recurrences implied by the
rip-grammar~\cite{rip:09}. Instead we build our theory as in
\cite{Li:10} around the concept of shapes, which we ``color'', in
order to rule out certain (bad) inflation scenarios. Passing from
shapes to refined shapes changes the shape-grammar as well as the
underlying generating functions. The refined shapes are key to the
generating functions since the collapsing of stems preserves vital
information of the interaction structure. It is therefore not
surprising that a shape induces joint structures via inflation, see
Theorem~\ref{T:RelateShape}.

\begin{figure}
\begin{center}
\includegraphics[width=0.7\textwidth]{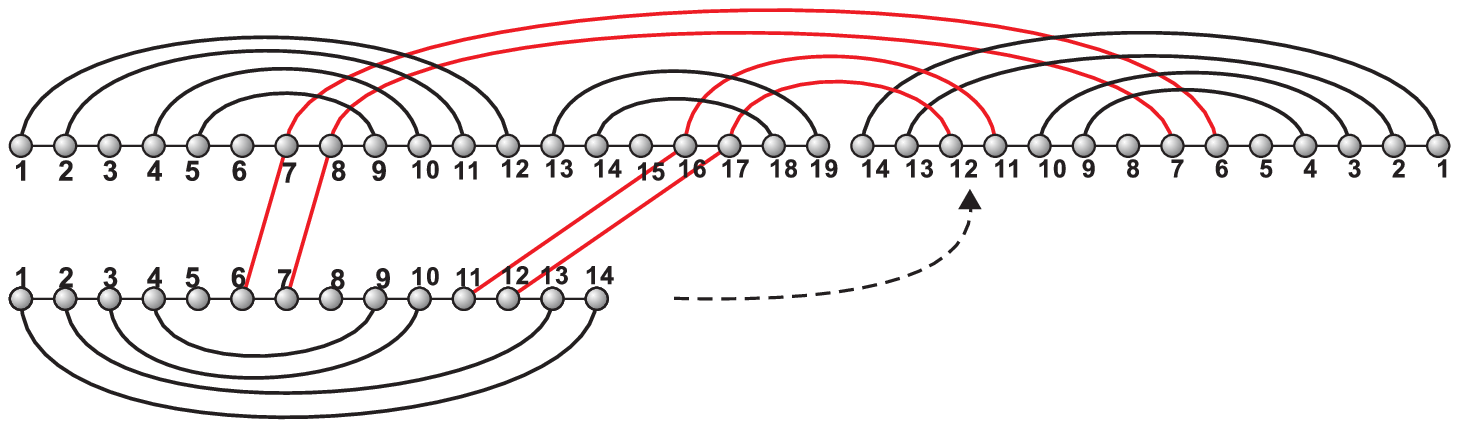}
\caption{Interpretation of joint structures as $3$-noncrossing structures.
A joint structure (left)
can be represented as $3$-noncrossing structure (top) by rotating the bottom
sequence around its endpoint
and aligning the two backbones and finally drawing all exterior arcs (red)
in the upper halfplane.}
\label{F:RNAiToCross}
\end{center}
\end{figure}

\begin{figure}
\begin{center}
\includegraphics[width=0.5\textwidth]{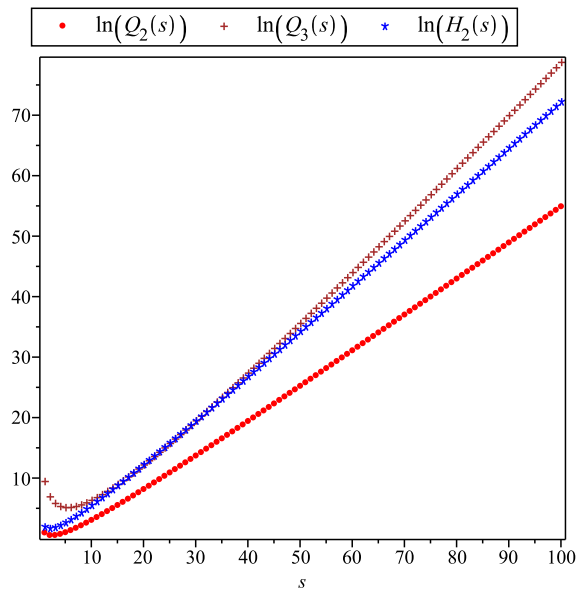}
\caption{How joint structures ``fit'' in:
we display the numbers of secondary structures (red), joint structures (blue) and
$3$-noncrossing pseudoknots structures (brown). All structure classes are canonical
and exhibit arc-length greater or equal to four.}
\label{F:compare}
\end{center}
\end{figure}

As canonical joint structures of arc-length at least four constitute
a novel combinatorial class it is of interest to compare them with
the classes of RNA secondary structures (having generating function
$\sum_nQ_2(s)z^s$) and $3$-noncrossing pseudoknots structures
($\sum_nQ_3(s)z^s$). Here a $3$-noncrossing structure has a diagram
representation in which there are no three mutually crossing arcs.
Indeed, clearly, RNA secondary structures are joint structures
without any exterior arcs. Furthermore any joint structure can be
interpreted as a particular $3$-noncrossing structure, by rotating
the bottom structure around its endpoint by 180 degrees, then
aligning the two backbones and drawing all exterior arcs in the
upper halfplane, see Fig.~\ref{F:RNAiToCross}.
For long sequences the numbers of canonical secondary structures,
$Q_2(s)$~\cite{Schuster}, joint structures $H_2(s)$ and
$3$-noncrossing pseudoknots structures
$Q_3(s)$, all having arc-length at least four~\cite{Modular} we find
\begin{eqnarray*}
Q_2(s) &\sim & 1.4848 \ s^{-\frac{3}{2}}\, 1.8489^s\\
H_2(s) &\sim & 3.5161 \ s^{-\frac{3}{2}}\, 2.1801^s \\
Q_3(s) &\sim & 5546 \ \ \ s^{-5}\ \, 2.5410^s,
\end{eqnarray*}
see also Fig.~\ref{F:compare}.

We can report that joint structures resemble features of secondary
structures as well as $3$-noncrossing
structures. Indeed, as it is the case for secondary structures, they can be
MFE-folded in polynomial time and as RNA pseudoknot structures they exhibit
crossing arcs and are truly shape-based
structure class. However, in contrast to $3$-noncrossing structures, refined shapes
have algebraic generating functions (as opposed to $D$-finite ones) and satisfy
simple recurrences.

Let us finally outline future research: with this paper the combinatorics of
joint structures is completed. The next step is to study their topology,
i.e.~understanding how joint structures filter via topological genera and
boundary components. This program means to pass from the deformation retracts studied
here to fat-graphs and their associated surfaces.


\section{Appendix}\label{S:facts}


\subsection{Singularity analysis}

In light of the fact that explicit formulas for the coefficients of a generating
function can be very complicated or even impossible to obtain, we estimate of the
coefficients in terms of the exponential factor and the subexponential factor.
Singularity analysis gives a framework that allows to extract the
asymptotics information of these coefficients. The key to obtain the
asymptotic formulas about the coefficients of a generating function is
its dominant singularities.
The theorem of Pringsheim~\cite{Flajolet:07a,Tichmarsh:39} guarantees
that a combinatorial generating function with nonnegative coefficients has
its radius of convergence as its dominant singularity. Furthermore for all
our generating functions it is the unique dominant singularity.
The derivation of exponential growth rates and subexponential factors from
singular expansions of generating functions mainly rely on the transfer
theorems~\cite{Flajolet:07a}.

To be precise, we say a function $f(z)$ is $\Delta_\rho$ analytic at
its dominant singularity $z=\rho$, if it analytic in some domain
$\Delta_\rho(\phi,r)=\{ z\mid \vert z\vert < r, z\neq \rho,\, \vert
{\rm Arg}(z-\rho)\vert >\phi\}$, for some $\phi,r$, where $r>|\rho|$
and $0<\phi<\frac{\pi}{2}$. We use the notation
\begin{equation*}
\left(f(z)=\Theta\left(g(z)\right) \
\text{\rm as $z\rightarrow \rho$}\right)\  \Longleftrightarrow \
\left(f(z)/g(z)\rightarrow c\ \text{\rm as $z\rightarrow \rho$}\right),
\end{equation*}
where $c$ is some constant. Let $[z^n]f(z)$ denote the coefficient
of $z^n$ in the power series expansion of $f(z)$ at $z=0$. Since the
Taylor coefficients have the property
\begin{equation*}
\forall \,\gamma\in\mathbb{C}\setminus 0;\quad [z^n]f(z)=\gamma^n
[z^n]f\left(\frac{z}{\gamma}\right),
\end{equation*}
We can, without loss of generality, reduce our analysis to the case
where $z=1$ is the unique dominant singularity. The following theorems
transfer the asymptotic expansion of a function around its unique
dominant singularity to the asymptotic of the function's
coefficients.
\begin{theorem}\label{T:transfer1}\cite{Flajolet:07a}
Let $f(z)$ be a $\Delta_1$ analytic function at its unique dominant
singularity $z=1$. Let
$$g(z)=(1-z)^{\alpha}\log^{\beta}\left(\frac{1}{1-z}\right)
,\quad\alpha,\beta\in \mathbb{R}.
$$
That is we have in the intersection of a neighborhood of $1$
\begin{equation}\label{E:transfer1}
f(z) = \Theta(g(z)) \quad \text{\it for } z\rightarrow 1.
\end{equation}
Then we have
\begin{equation}
[z^n]f(z)= \Theta\left([z^n]g(z)\right).
\end{equation}
\end{theorem}
\begin{theorem}\label{T:transfer2}\cite{Flajolet:07a}
Suppose $f(z)=(1-z)^{-\alpha}$, $\alpha\in\mathbb{C}\setminus
\mathbb{Z}_{\le 0}$, then
\begin{equation}
\begin{aligned}
[z^n]\, f(z)  \sim & \frac{n^{\alpha-1}}{\Gamma(\alpha)}\left[
1+\frac{\alpha(\alpha-1)}{2n}+\frac{\alpha(\alpha-1)(\alpha-2)(3\alpha-1)}
{24 n^2}+\right. \\
&\qquad  \quad
\left.\frac{\alpha^2(\alpha-1)^2(\alpha-2)(\alpha-3)}{48n^3}+
O\left(\frac{1}{n^4}\right)\right].
\end{aligned}
\end{equation}
\end{theorem}


\subsection{Symbolic Enumeration}

Symbolic enumeration~\cite{Flajolet:07a} plays an important role in the
computations of generating functions. We first introduce the notion of a
combinatorial class. Let
$\mathbf{z}=(z_1,\ldots,z_d)$ be a vector of $d$ formal variables
and $\mathbf{k}=(k_1,\ldots,k_d)$ be a vector of integers of the
same dimension. We use the simplified notation
\begin{equation*}
\mathbf{z}^{\mathbf{k}}\colon=z_1^{k_1} \cdots z_d^{k_d}.
\end{equation*}
\begin{definition}
A combinatorial class of $d$ dimension, or simply a class, is an
ordered pair $(\mathcal{A},w_{\mathcal{A}})$ where $\mathcal{A}$ is
a finite or denumerable set and a size-function
$w_{\mathcal{A}}\colon \mathcal{A}\longrightarrow \mathbb{Z}_{\geq
0}^d$ satisfies that $w_{\mathcal{A}}^{-1}(\mathbf{n})$ is finite
for any $\mathbf{n}\in\mathbb{Z}_{\geq 0}^d$.
\end{definition}
Given a class $(\mathcal{A},w_{\mathcal{A}})$, the size of an
element $a\in\mathcal{A}$ is denoted by $w_{\mathcal{A}}(a)$, or
simply $w(a)$. We consistently denote by $\mathcal{A}_{\mathbf{n}}$
the set of elements in $\mathcal{A}$ that have size $\mathbf{n}$ and
use the same group of letters for the cardinality
$A_{\mathbf{n}}=\vert \mathcal{A}_{\mathbf{n}}\vert$. The sequence
$\{A_{\mathbf{n}}\}$ is called the counting sequence of class
$\mathcal{A}$. The generating function of a class
$(\mathcal{A},w_{\mathcal{A}})$ is given by
\[
\mathbf{A}(\mathbf{z})=\sum_{a\in\mathcal{A}}\mathbf{z}^{w_{\mathcal{A}}(a)}=
\sum_{\mathbf{n}}A_{\mathbf{n}}\, \mathbf{z}^{\mathbf{n}}.
\]
There are two special classes: $\mathcal{E}$ and $\mathcal{Z}_i$
which contain only one element of size $\bf{0}$ and $\mathbf{e}_i$,
respectively. In particular, the generating functions of the classes
$\mathcal{E}$ and $\mathcal{Z}_i$ are
\[
{\bf{E}}(\mathbf{z})=1\quad\text{and}\quad {\bf{Z}}_i(\mathbf{z})=z_i.
\]
Next we introduce some basic constructions that constitute the core of a
specification language for combinatorial structures. Let $\mathcal{A}$ and
$\mathcal{B}$ be combinatorial classes of $d$ dimension. Suppose
$\mathcal{A}_{i}$ are combinatorial classes of $1$ dimension. We define
\begin{itemize}
\item $(\mathcal{A}_1,\mathcal{A}_2):=\{c=(a_1,a_2)\mid
a_i\in\mathcal{A}_i\}$ and for
$c=(a_1,a_2)\in(\mathcal{A}_1,\mathcal{A}_2)$
\[
w_{(\mathcal{A}_1,\mathcal{A}_2)}(c)=
(w_{\mathcal{A}_1}(a_1),w_{\mathcal{A}_2}(a_2))),
\]
\item $\mathcal{A}+\mathcal{B}:=\mathcal{A}\cup\mathcal{B}$, if
$\mathcal{A}\cap\mathcal{B}=\varnothing$ and for
$c\in\mathcal{A}+\mathcal{B}$,
\[
w_{\mathcal{A}+\mathcal{B}}(c)=\left\{
\begin{aligned}
&w_{\mathcal{A}}(c)\quad \text{if}\ c\in\mathcal{A}\\
&w_{\mathcal{B}}(c)\quad \text{if}\ c\in\mathcal{B},
\end{aligned} \right.
\]
\item $\mathcal{A}\times\mathcal{B}:=\{c=(a,b)\mid
a\in\mathcal{A},b\in\mathcal{B}\}$ and for
$c\in\mathcal{A}\times\mathcal{B}$,
\[
w_{\mathcal{A}\times\mathcal{B}}(c)=
w_{\mathcal{A}}(a)+w_{\mathcal{B}}(b),
\]
\item
$\textsc{Seq}(\mathcal{A}):={\mathcal{E}}+\mathcal{A}+(\mathcal{A}\times\mathcal{A})+
(\mathcal{A}\times\mathcal{A}\times\mathcal{A})+\cdots$.
\end{itemize}
Plainly, $\textsc{Seq}(\mathcal{A})$ defines a proper combinatorial
class if and only if $\mathcal{A}$ contains no element of size $0$.
We immediately observe
\begin{proposition}\label{P:Symbolic}
Suppose $\mathcal{A}$, $\mathcal{B}$ and $\mathcal{C}$ are
combinatorial classes of $d$ dimension having the generating functions
$\mathbf{A}(\mathbf{z})$, $\mathbf{B}(\mathbf{z})$ and
$\mathbf{C}(\mathbf{z})$. Let $\mathcal{A}_{i}$ be combinatorial
classes of $1$ dimension having the generating functions
$\mathbf{A}_i(z)$. Then\\
{\sf(a)} $\mathcal{C}=(\mathcal{A}_1,\mathcal{A}_2,\ldots,\mathcal{A}_d)
\Longrightarrow \mathbf{C}(\mathbf{z})=\mathbf{A}_1(z_1)\,
\mathbf{A}_2(z_2)\ldots \mathbf{A}_d(z_d)$\\
{\sf (b)} $\mathcal{C} =\mathcal{A}+\mathcal{B} \Longrightarrow
\mathbf{C}(\mathbf{z})=\mathbf{A}(\mathbf{z})+\mathbf{B}(\mathbf{z})$\\
{\sf (c)} $\mathcal{C}=\mathcal{A}\times\mathcal{B}\Longrightarrow
\mathbf{C}(\mathbf{z})=\mathbf{A}(\mathbf{z})\cdot\mathbf{B}(\mathbf{z})$\\
{\sf (d)} $\mathcal{C}=\textsc{Seq}(\mathcal{A})\Longrightarrow
\mathbf{C}(\mathbf{z})=\frac{1}{1-\mathbf{A}(\mathbf{z})}.$
\end{proposition}

\bibliographystyle{elsarticle-num}

\end{document}